# Solving Groundwater Flow Equations Using Gradually Varied Functions

Li Chen
Department of Computer Science and Information Technology
University of the District of Columbia
Email: *lchen@udc.edu*

Xun-Hong Chen
School of Natural Resources
University of Nebraska-Lincoln
Email: *xchen2@unl.edu*

**Abstract**
Finite difference method and finite element method are popular methods for solving groundwater flow equations. This paper presents a new method that uses gradually varied functions to solve such equation. In this paper, we have established a mathematical model based on gradually varied functions for groundwater data volume reconstruction. These functions do not rely on the rectangular Cartesian coordinate system. A gradually varied function can be defined in a general graph or network. Gradually varied functions are suitable for arbitrarily shaped aquifers. Two types of models are designed and implemented for real data processing: (1) the gradually varied model for individual (time) groundwater flow data, (2) the gradually varied model for sequential (time) groundwater flow data. In application, we used two sets of real data and one set of experimental data to test our methods.

## 1. Introduction

The current method of flow estimation mainly uses the ground flow equation, which is a partial differential equation. Computer source codes such as MODFLOW solves 2D equations and pass the data vertically to form a 3D volume. The research on the 3D models of groundwater flow has fundamental and practical importance to hydrogeology. A method used to establish a true 3D Groundwater flow will be very useful to groundwater.

Much research has already been done to find a discrete model for the groundwater flow equations. The research is mainly based on numerical methods and analytical methods. The finite difference method (FDM) and the finite element method (FEM) are popular in this area. Pruist et al. [22] has indicated that FEM has advantages on local refinement of grid (adaptive mesh generation) due to non-rectangular grids, good accuracy, stability, representation of the spatial variation of anisotropy. However, its computational cost is much larger and relative less manageable in application. In fact, groundwater industry is not like the automobile industry, there is no much need for a good-looking smoothed groundwater level surface.

On the other hand, FDM has simplicity of theory and algorithm, plus easiness of application; however it has the problems of inefficient refinement of grid and poor geometry representation. This is because of the strict use of rectangular grids. In addition, there is no standard method to implement the Neumann boundary condition.

The gradually varied function is supposed to pick up the advantages and overcome the disadvantages in FEM and FDM. So our first task is to investigate the suitableness of gradually varied functions for groundwater data. Then, we must find a connection between the flow equations and gradually varied functions. We also need to design an input data format to store the data in a database.

A gradually varied function is for the discrete system where a high level of smoothness is not a dominant factor. It can be used in any type of decomposition of the domain. It is more flexiable than rectangle-cells used in MODFLOW and triangle-cells used in FEFLOW. Because gradual variation does not have strict system requirements, the other mathematical methods and the artificial intelligence methods can be easily incorprated into this method to seek a better solution. Based on the boundary conditions or constraints of the groundwater aquifer, the constraints could be in explicit forms or in differential forms such as the diffusion equations. The gradually varied function exists based on the following theorem: The necessary and sufficient condition of the existance of a gradually varied function is that the change of values in any pair of sample points is smaller or equal to the distance of the points in the pair.

In order to use discrete methods for groundwater modeling, we have completed the following tasks: (1) Groundwater flow equations and its discrete forms, (2) Gradually varied functions for groundwater data, (3) Real data preparation, (4) Algorithms, especially fast design, and (5) Real data processing and applications.

## 2. Background and Related Research

The research of groundwater flow is one of the major topics in subsurface hydrology [1] [2] (this sentence does not make sense and just repeats itself). The groundwater flow equation based on Darcy's Law usually describes the movement of groundwater in a porous medium such as aquifers. It is known in mathematics as the diffusion equation. The Laplace equation (for a steady-state flow) is a simplification of the diffusion equation under the condition that the aquifer has a recharging and/or discharging boundary. The conservation of mass states that for a given increment of time ($\Delta t$) "the difference between the mass flowing in across the boundaries, the mass flowing out across the boundaries, and the sources within the volume, is the change in storage."[2]

$$\frac{\Delta M_{stor}}{\Delta t} = \frac{M_{in}}{\Delta t} - \frac{M_{out}}{\Delta t} - \frac{M_{gen}}{\Delta t} \qquad (1)$$

Its differential form is

$$\frac{\partial h}{\partial t} = \alpha \left[ \frac{\partial^2 h}{\partial x^2} + \frac{\partial^2 h}{\partial y^2} + \frac{\partial^2 h}{\partial z^2} \right] - G \qquad (2)$$

Where $h$ is the hydraulic head, $\alpha$ is a mathematical term and it is the function of aquifer storage coefficient $S$, aquifer thickness $b$, and aquifer hydraulic conductivity $K$, and G is the sink/source term such as groundwater pumping, recharge etc. To solve equation (2), a grid method is usually used, such as the finite difference or finite element method [3] [4]. Other methods including the analytic element method attempt to solve the equation exactly, but need approximations of the boundary conditions [5][6]; they are mainly used in academic and research labs.

In practice, numerical methods have been commonly used to solve these equations. The most popular system is called MODFLOW, which was developed by USGS [7]. MODFLOW is based on the finite-difference method on rectangular Cartesian coordinates. MODFLOW can be viewed as a "quasi 3D" simulation since it only deals with the vertical average (no $z$-direction derivative). Flow calculations between the 2D horizontal layers use the concept of leakage.

Gradual variation is a discrete method that can build on any decomposition or networking. It was originally introduced to solve image processing problems and discrete surface reconstruction [8][9][21]. See Fig 1. When a boundary is known, it can be used to fit (solve) for the interal points in any type of linked connection.

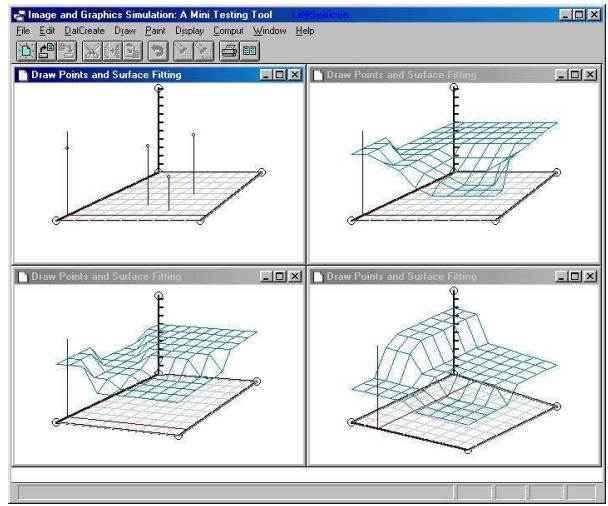

Fig. 1 Examples of gradually varied functions

The original idea of using gradually varied surfaces in this project is to use Darcy's Law or differential constrains to determine the hydraulic head value of the unknown points instead of the random selection of the construction of gradually varied surfaces when there are more than one possible selections [10]. When the determinations of the values are uncertain, one can try to use the artificial intelligence methods such as neural networks and genetic algorithms that help us to find near optimal solutions. These types of studies have already been done by many researchers in subsurface hydrology field [11][12].

Gradual variation is a discrete method that can be built on any graph. The gradually varied surface is a special discrete surface. We now introduce this concept.

The Concept of Gradual Variation: Let function f: D$\rightarrow${1, 2,…,n}, if a and b are adjacent in D implies |f(a)- f(b)| $\leq$1, point (a,f(a)) and (b,f(b)) are said to be gradually varied. A 2D function (surface) is said to be gradually varied if every adjacent pair are gradually varied.

Discrete Surface Fitting: Given J$\subseteq$D, and f: J$\rightarrow${1,2,…n} decide if there is a F: D$\rightarrow${1,2,…,n} such that F is gradually varied where f(x)=F(x), x in J.

**Theorem 1** (Chen, 1989) The necessary and sufficient conditions for the existence of a gradually varied extension F is: for all x,y in J, d(x,y)$\geq$ |f(x)-f(y)|, where d is the distance between x and y in D.

The above theorem can be used for a single surface fitting if the condition in the theorem is satisfied. The problem is that the sample data does not satisfy the condition of fitting. So the original algorithm cannot be used directly for individual surface fitting.

On the other hand, current research still showed considerable interest in establishing new modeling methods for groundwater flow [13-16]. Both the US and DC governments are concerned with the future of groundwater flow research [17][18]. This project on the 3D models of groundwater flow has fundamental and practical importance to hydrogeology. A method used to establish a true 3D model for groundwater flow will be very useful to groundwater related research in DC and in other urban areas.

### 3. Ground Water Surface Fitting

We have done extensive research on data collection and initial data reconstruction using gradually varied functions that are for the discrete system where a high level of smoothness is not a dominant factor. Since they can be used on any type of decomposition of the domain, gradually varied functions are more flexible than rectangle-cells used in MODFLOW and triangle-cells used in FEFLOW. Because gradual variation does not have strict system requirements, the other mathematical methods and the artificial intelligence methods can be easily incorporated into this method to seek a better solution.  Based on the boundary conditions or constraints of the groundwater aquifer, the constraints could be in explicit forms or in differential forms such as the diffusion equations. Again we select the gradually varied function because of  the following theorem: The necessary and sufficient condition of the existence of a gradually varied function is that the change of values in any pair of sample points is smaller or equal to the distance of the points in the pair.

### 3.1 Individual surface fitting

Directly use Theorem 1 to fit the data will need to evenly divide the grade level. An algorithm (Algorithm A) based on the sample point contribution to the fitting point is created. The core part of the program is given in the following code. This algorithm can produce the surface but do not give good fit in some cases. We have showed the result in the previous report.

```
Algorithm A:
•       for (k=0;k<nGuildPoints;k++){
•               ii=(latIndex[k]-latMin)/latDet;
•               jj=(longtIndex[k]-longtMin)/longtDet;

•               distance=sqrt((ii -i)*(ii-i)+(jj -j)*(jj-j));
•                 temp_j=abs((array[i][j] - dat[k][time]))/Ratio-distance;
                        if(temp_j>0){ // not satisfy gvs condition

•               if( array[i][j] > dat[k][time])
•                       temp=-temp_j *Ratio ;
•               else
•                       temp= temp_j *Ratio;
•               array[i][j]=array[i][j]+temp ;
•       }
•
```

This algorithm will work fine when the data values do not jump largely for closer neighborhood.

A new algorithm was developed recently to overcome these problems [23]. In [23], a systematic digital-discrete method for obtaining continuous functions with smoothness to a certain order ($C^n$) from sample data is designed. This method is also based on gradually varied functions.  The new algorithm tries to search for a best solution of the fitting. We have added the component of the classical finite difference method. The result of the new method will be shown in following sections. This method is independent from existing popular methods such as the cubic spline method and the finite element method. The new digital-discrete method has considerable advantages for a large amount of real data applications. This digital method also differs from other classical discrete method that usually uses triangulations. This method can potentially be used to obtain smooth functions such as polynomials through its derivatives $f^{(k)}$ and the solution for partial differential equations such as harmonic and other important equations.

After different derivatives are obtained, we can use Taylor expansion to update the value of the gradually varied fitted function (at $C^{(0)}$). In fact, at any order $C^{(k)}$, we can update it using a higher order of derivatives.

The Taylor expansion is based on the formula of the Taylor series, which has the following generalized form:

$$f(x_1, \cdots, x_d) = \sum_{n_1=0}^{\infty} \cdots \sum_{n_d=0}^{\infty} \frac{(x_1 - a_1)^{n_1} \cdots (x_d - a_d)^{n_d}}{n_1! \cdots n_d!} \left( \frac{\partial^{n_1 + \cdots + n_d} f}{\partial x_1^{n_1} \cdots \partial x_d^{n_d}} \right)(a_1, \ldots, a_d). \quad (3)$$

For example, for a function that depends on two variables, $x$ and $y$, the Taylor series of the second order using the guiding point $(x_0, y_0)$ is:

$$f(x, y) \approx f(x_0, y_0) + (x - x_0) f_x(x_0, y_0) - (y - y_0) f_y(x_0, y_0) \quad (4)$$

There are several ways of implementing this formula. We have chosen the ratio of less than half of the change. An iteration process is designed to make the new function converge.

**Algorithm 2.**
The new algorithm tries to search for the best fit. We have added a component of the classical finite difference method. Start with a particular dataset consisting of guiding points defined on a grid space. The major steps of the new algorithm are as follows (This is for 2D functions. For 3D function, we would only need to add a dimension):

    **Step 1:** *Load guiding points. In this step we load the data points with observation values.*
    **Step 2:** *Determine the resolution. Locate the points in grid space.*
    **Step 3**: *Function extension according to Theorem 2.1. This way, we obtain gradually varied or near gradually varied (continuous) functions. In this step the local Lipschitz condition is used.*
    **Step 4:** *Use the finite difference method to calculate partial derivatives. Then obtain the smoothed function.*
    **Step 5:** *Some multilevel and multi resolution method may be used to do the fitting when data set is large.*

## 3.2 Sequential surface fitting and involvement of the flow equation

The individual surface fitting in the above code is not an adequate method since the relationship of groundwater surfaces at each time is not accounted for in the calculation. In particular, we have to use the flow equation in the sequential surface calculation.

According to flow equation (2),
$$h2-h1 = alpha*(h2(x-1,y)+h2(x+1,y)-2h2(x,y)+ h2(x,y-1)+h2(x,y+1)-2h2(x,y)) -G \quad (5)$$

We can let
$$f4= (h2(x,y)-h1(x,y)+G)/alpha + 4*h2(x,y) = h2(x-1,y)+h2(x+1,y)+ h2(x,y-1)+h2(x,y+1). \quad (6)$$

f4 can also be viewed as the average of 4 times the h values at the center point. Using the gradually varied function, the average is ((h2-h1+G)/alpha + 4*h2(x,y))/4. Based on the properties of gradually varied functions, the maximum difference is three. For example, h2_old(x-1),y) is bigger , then h2_new(x-1),y) should be smaller. The iterating procedure is to update the center point after the first fit.

## 4. Experiments and Real Data Processing

We have applied the above algorithms to three sets of problems. The first set is the data from a region in Northern Virginia.
The second set is experimental data around a pumping well. The third set of the data is also pumping well but it is real data.

### 4.1 Lab Data Reconstruction

In this data set, we have simulated a pumping well. We select 40 observation locations surrounding the well to get the log data. Those points are selected along with eight transects surrounding the pumping well. Five observation points are chosen in each transect.

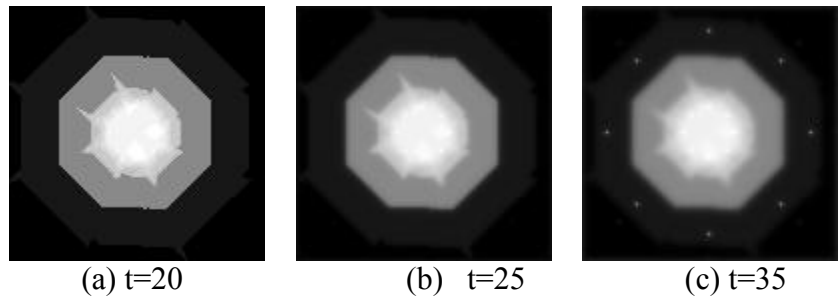

(a) t=20    (b) t=25    (c) t=35

Fig. 3. Lab data experimentation. (a) t=20  (b)  t=25           (c) t=35. there are eight dots remaining on the picture. That means we need to select more   suitable  *G* value.

This example shows that our algorithm worked very well on data reconstruction. The correctness of our algorithm is partially verified. If 16 directions with 80 point samples are chosen, the resulting picture will be much circular.

**4.2 Real Data reconstruction based on one pumping well**

This example contains on 9 samples of real data points for a pumping well.

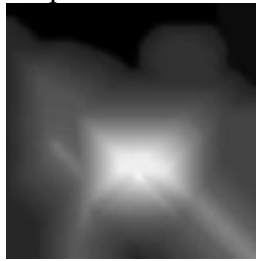

Fig. 4. A pumping well with 9 sample points.

The data distribution has a rectangular pattern. It is caused by the reconstruction process that was made based on rectangle cells. When sample points are small amount and they are not selected randomly, the data will appear to be rectangle shape.

4.3 Real data reconstruction for a large region
Fig. 5 shows the results that were produced using an individualized-fit algorithm to fit the initial data set. This algorithm is also made by the rough graduate varied surface fitting by scanning through the fitting array. There are many clear boundary lines in the images.

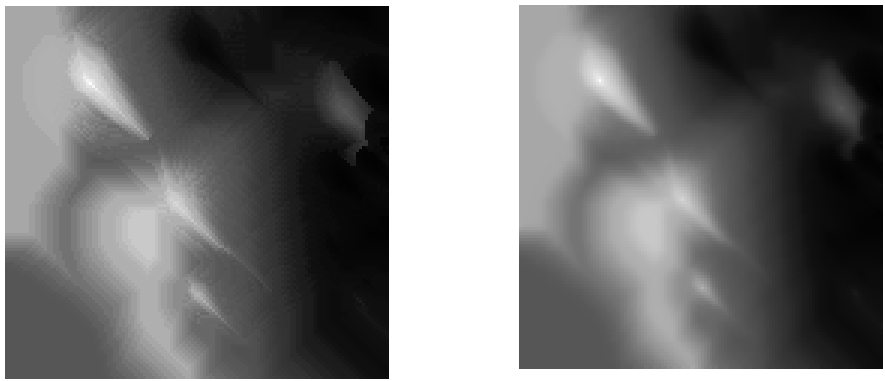

Fig. 5 Gradually varied functions for initial surface

Fig.6 shows the sequential surface fitting and the water flow equation updating. Individual surface fitting results are shown in (a), (b), and (c). Starting with fitted surface at each time the process will be faster to converge. It will not affect the final result if there are enough iterations.

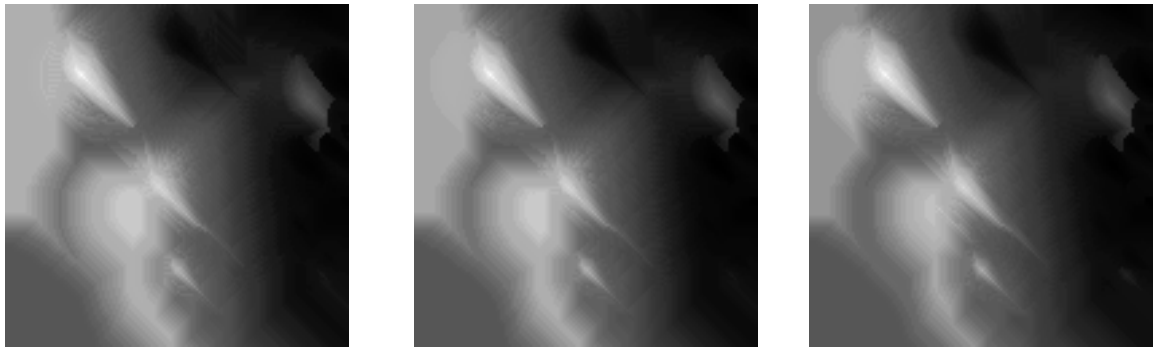

(a) day1  (b) day 30  (c) day 50

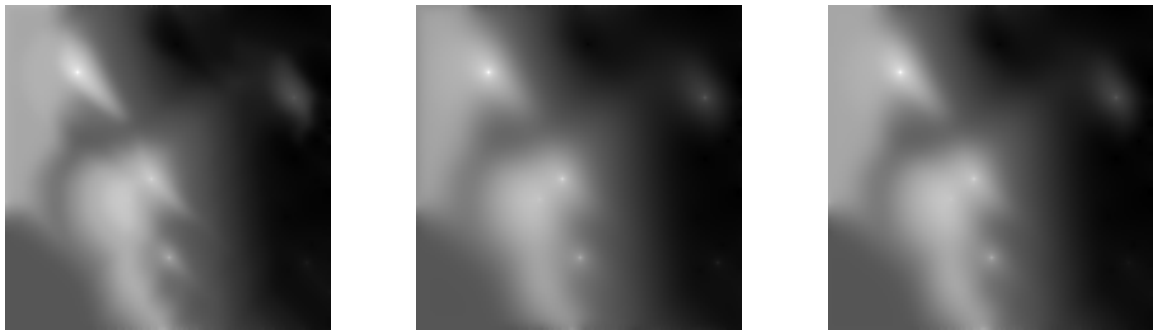

(d) Day 3 with flow equation iteration   (f) Day 30 with iteration   (g) Day 50 with iteration

Fig 6. VA Groundwater distribution calculated by gradually varied surfaces and the flow equation.

To test the correctness of the result, we have located the data points in the rectangle area in GoogleMap and Find latitude and Longitude. http://findlatitudeandlongitude.com/
We also tested the following fitting points that show the corrected correspondences. Data preparation is a very important part of this project. Please see appendix or [20] for more information. A web programming language PHP was used to build a web application to access groundwater log data in Virginia and Maryland. Data is stored in MySQL databases.

Selected Points used in reconstruction
One can find the location at

| | | |
|---|---|---|
| 4.65 | 36.62074879 | -76.10938540 |
| 75.37 | 36.92515020 | -77.17746768 |
| 6.00 | 36.69104276 | -76.00948530 |
| 175.80 | 36.78431615 | -76.64328700 |
| 168.33 | 36.80403855 | -76.73495750 |
| 157.71 | 36.85931567 | -76.58634110 |
| 208.26 | 36.68320624 | -76.91329390 |
| 7.26 | 36.78737704 | -76.05153760 |

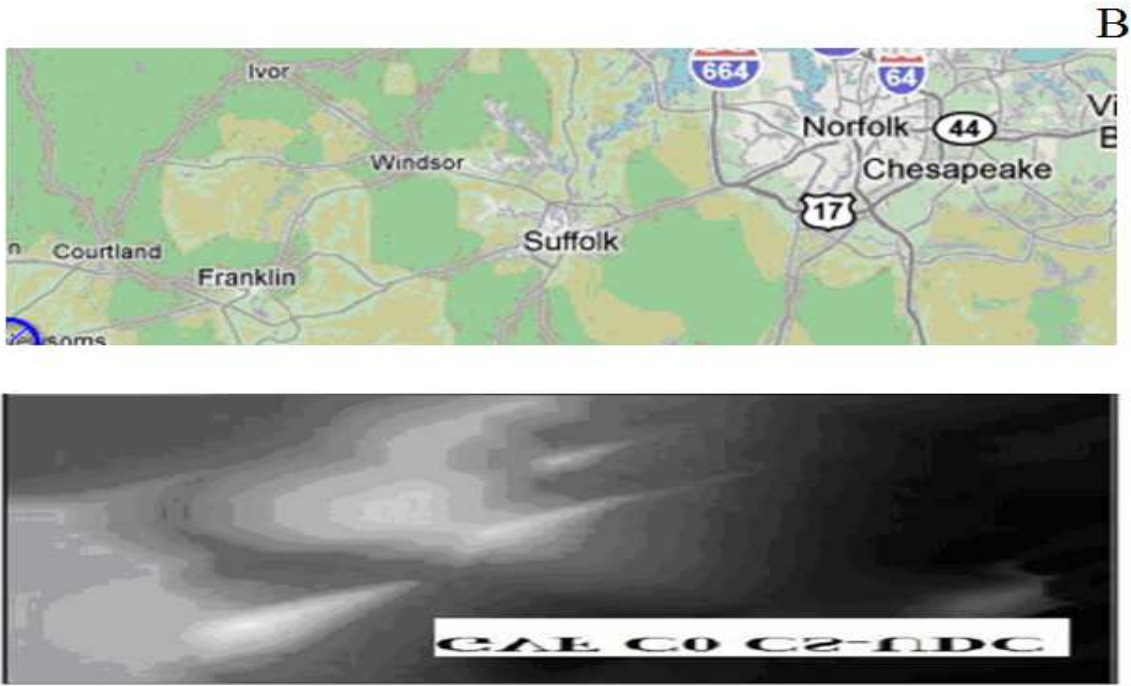

Dimensions (Latitude, Longitude)
A = (36.62074879, -77.17746768)
B = (36.92515020, -76.00948530)

Fig. 7 The map and ground water data

Fig. 7 shows a good match between the ground water data and the region geographical map. The brighter pixels mean the deeper distance from the surface. In mountain area, the groundwater level is lower in general. Some mismatches may be caused by not having enough sample data points (observation wells).

**Future work:** to ensure the accuracy of the calculation, we will add the finite element method to our research. We also want to use MODEFLOW to calculate local and small-region flow, and to use gradual variation to compute regional or global data.

*Acknowledgement:* This research was supported by USGS seed grants. The author expresses thanks to Dr. William Hare and members of UDC DC Water Resources Research Institute for their help. The author also thank Mr. Travis L. Branham for data collection. This paper is based on the report entitled: Li Chen, Gradual Variation Analysis for Groundwater Flow of DC, DC Water Resources Research Institute Final Report 2009 [24].

.